# Modified FEA and ExtraTree algorithm for transformer Green's function modelling

Xuhao Du[1]; Jie Pan[1]

[1] The University of Western Australia, School of Mechanical and Chemical Engineering, Perth, Australia

**ABSTRACT**

The Green's function of a transformer is essential for prediction of its vibration. As the Green's function cannot be measured directly and completely, the finite element analysis (FEA) is typically used for its estimation. However, because of the complexity of the transformer structure, the calculations involved in FEA are time consuming. Therefore, in this paper, a method based on FEA modified by an extremely random tree algorithm call ExtraTree is proposed to efficiently estimate the Green's function of a transformer. A subset of the frequency response functions from FEA will be selected by a genetic algorithm that can well present the structural variation. The FEA calculation time can be reduced by simply calculating the frequency response functions on this subset and predicting remainder using the trained ExtraTree model. The errors introduced in this method can be estimated from the corresponding frequency and the genetic algorithm error.

Keywords: Green's function, Finite element analysis, Optimisation, ExtraTree, Transformer
Subjects Number(s): 51.4

## 1. INTRODUCTION

Calculation of the mechanical properties of a power transformer, such as the Green's function, has been a challenge to researchers for a few decades. The Green's function is the assembly of all the frequency response functions (FRF) of the structure. It can be used for vibration prediction by combining the force density inside the core. However, because of the complicated friction model, ensemble type, and boundary conditions, there is still no analytical model that can represent the Green's function of a transformer structure.

In the past few years, much effort has been put into calculating the FRFs of a transformer core. In 2007, Duhamel proposed using finite element analysis (FEA) to compute the Green's function of any linear homogeneous medium, relying on wave propagation theory in periodic media. This method was verified by calculating the Green's function of a periodic array and acoustics element, whose Green's functions are already known (1). In 2011, Pan *et al.* used FEA to calculate the first few Eigen modes of an actual transformer and compared them to experimental results (2). They took the transformer core as a complete solid structure and calculated the fundamental frequency with 8% deviation compared to the experimental result.

On the other hand, as transformers are typical lamination structures, their lamination properties have been studied for several decades. In 1964, Walker *et al.* studied the minimum pressure for tightening the lamination core (3). In 2004, Masti *et al.* reported that part of a transformer core's vibration is due to the relative motion of the lamination, which has a strong influence on a structure's resonant behaviour. By comparing the natural frequencies of the lamination structure and core as a complete body, they used FEA to simulated the core and found that the natural frequencies displayed a huge reduction, and this was verified by the experimental results (4).

As mentioned above, FEA is widely used in the Green's function modelling of power transformers. However, there is still error and huge consumption of both time and memory. Therefore, instead of a purely FEA calculation, a novel approach for modelling Green's functions of transformer is pro- posed by using an ExtraTree algorithm to train a subset of the FEA-FRFs. The ExtraTree algorithm is a machine learning algorithm which is basically a bagging model of multiple decision trees (5). It has been demonstrated to be useful in reducing both the model's variation and bias on classification and regression problems. Because of the continuity and high damping ratio of the transformer structure, the response at a certain node may vary slightly compared to its neighbour, especially in the low frequency range, where the wavelength is comparable to the dimensions of the structure. Thus, ExtraTree is suitable for predicting the FRFs of a transformer.

---

[1] xuhao.du@research.uwa.edu.au

Before training the ExtraTree model, a subset of FRFs that can describe the complete structural variation need to be selected. In this work, this subset is obtained by using a genetic algorithm (GA). A GA is a mimic algorithm for optimisation, especially for a spatial arrangement. Because of its robustness and high efficiency, GAs have been recognised as a useful optimisation tool in active vibration and noise control. In 1999, Han & Lee placed his piezoelectric sensors and actuators with GA optimisation for active control of a cantilevered composite plate. The performance was significant with 15 dB reduction of the plate vibration. Their close loop tests also demonstrated the GA's robustness (6). In 2004, Li *et al.* used a GA to optimise the location of a PZT actuator in an active structural acoustics control system, and this resulted in a significant reduction of the sound emission in 2004 (7).

In this paper, a novel method is proposed for efficiently estimating the Green's function based on the FEA, the ExtraTree model, and a GA. A domain zone theory is proposed to explain the mechanism behind this method. Details will be introduced in the following sections.

## 2. METHODOLOGY

The procedure for this method is presented in Figure 1. To establish the FEA model, a preliminary measurement will be adopted as a reference for the modification of FEA parameters. Next, the full set of FRFs at a certain frequency, named the seed frequency, is calculated as the base dataset.

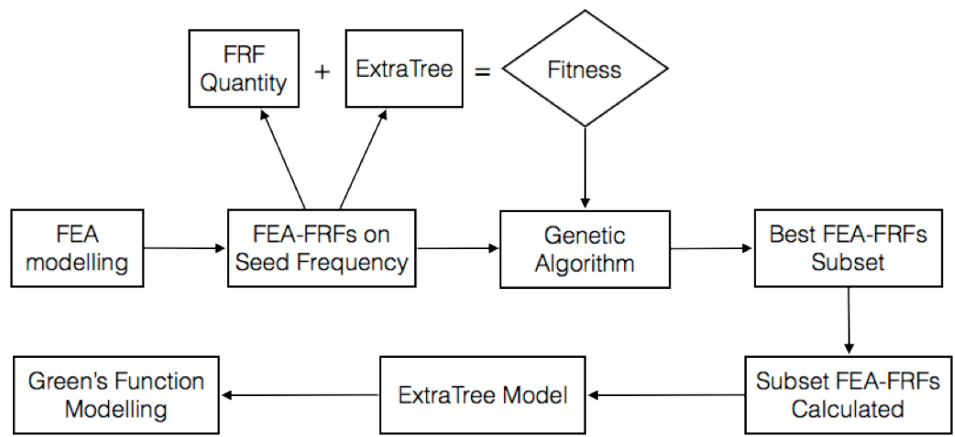

Figure 1 – A flowchart for the procedure

After obtaining the base dataset of the FEA-FRFs at the seed frequency, the GA will be applied to it. GAs are widely used in position optimisation and can be used to determine the optimal subset of FEA-FRFs for training the ExtraTree model. A GA is a mimic algorithm that simulates the evolutionary processes by selection, crossover and mutation. In this process, a criterion will be applied to each individual and its score under this criterion, called the fitness, will be calculated. As the number of steps of evolution increases, the fitness of the individual will improve. In this case, the fitness criterion is the combination of the FEA-FRF quantity of the selected subset and the performance of the ExtraTree models trained on the corresponding subset. After a few rounds of evolution, there will be one individual in the last generation that has the best fitness score, and this is the best GA subset for modelling the Green's function. The following is the pseudocode for the GA (8).

---------------------------------------------------Genetic Algorithm-----------------------------------------------

1. Generate the individuals with a certain population $P$ as the $1^{st}$ generation.
2. For $g = 1$ to $E$ (where $E$ is the set evolution time):
   (a) Input $g^{th}$ generation:
      i. Sort the fitness of the $g^{th}$ generation by combining the FEA-FRF used quantity and the ExtraTree performance.
      ii. Select the first $Pr$ individuals in the generation as the survivors, where $r$ is the ratio of those remaining.
      iii. Select the $Pr_s$ individuals in the generation as the survivors, where $r_s$ is the ratio of those randomly selected.
      iv. For $c = 1$ to $P*(1-r_s-r)$:
         A. Randomly select two survivors as the parents
         B. Crossover the parents to generate new individual $g_c$
      v. Run through the $g^{th}$ generations with $m$ percentage to mutate them.
      vi. The $g^{th}$ generation are generated completely.
3. Sort the fitnesses of the $E^{th}$ generation. The first individual is the best individual that is searched for.

-----------------------------------------------------------------------------------------------------------------

As a result, the best subset will be obtained, accompanied by an acceptable ExtraTree model. For obtaining the FRFs at other frequencies, the FRFs in the selected subset need to be calculated via FEA. Then the FEA result is trained using the ExtraTree algorithm to develop a Green's function model at the corresponding frequency. The values of the remainder of the base dataset can be estimated using the trained ExtraTree model with controlled error. ExtraTree is a binary tree-based machine learning algorithm which will be used in GA subset selection and Green's function modelling. Its pseudocode is as follows (8).

-------------------------------------------Extremely Random Tree Algorithm-----------------------------------------

1. For $b = 1$ to $B$:
   a. Draw a bootstrap sample Z of size N from the training data
   b. Grow a random-forest tree $T_b$ to the bootstrapped data, by recursively repeating the following steps for each terminal node of the tree, until the minimum node size $n_{min}$ is reached.
      i. Randomly select $m$ variables from the $p$ variables.
      ii. Randomly select the best variable and split-point among the $m$ variables.
      iii. Split the node into two daughter nodes.

2. Output the ensemble of trees $\{T_b\}_1^B$

$$f_{et}^B(x) = \frac{1}{B}\sum_{b=1}^{B} T_b(x) \quad (1)$$

-----------------------------------------------------------------------------------------------------------------

The details for conducting each step and algorithm will be presented in the following sections.

## 3. EXPERIMENT AND FEA

To establish an accurate FEA model, it is important to use a preliminary experiment as the reference for the model. In this section, the experiment setup and measurement results will be introduced and the result of FEA simulation will be presented.

### 3.1 Preliminary Experiment

The experiment was conducted on a single-phase transformer with fixed ground boundary conditions, as presented in Figure 2. The dimensions of this transformer are 500 mm in height, 250 mm in width and 50 mm in depth. The transformer is assembled using the multiple-layers method at the corner. There are two parallel windings and each one has 415 turns.

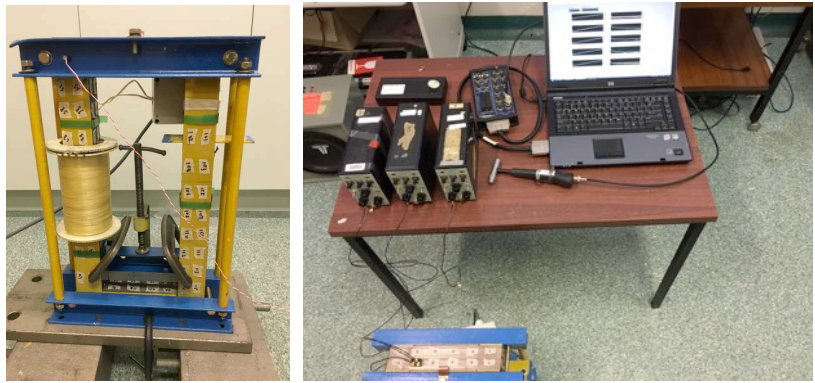

Figure 2 – The single-phase transformer and the data acquisition system

A B&K impact hammer was used to excite the transformer's vibration and its responses were picked up by a triaxial accelerometer and collected by an NI card with LabVIEW software. The frequency responses at multiple nodes were measured and the average frequency response of the transformer is shown in Figure 3. A few clear peaks can be observed from Figure 3 as the natural frequencies, which are listed in Table 1. The mode shapes of the first four natural frequencies were measured and the first on is presented in Figure 4.

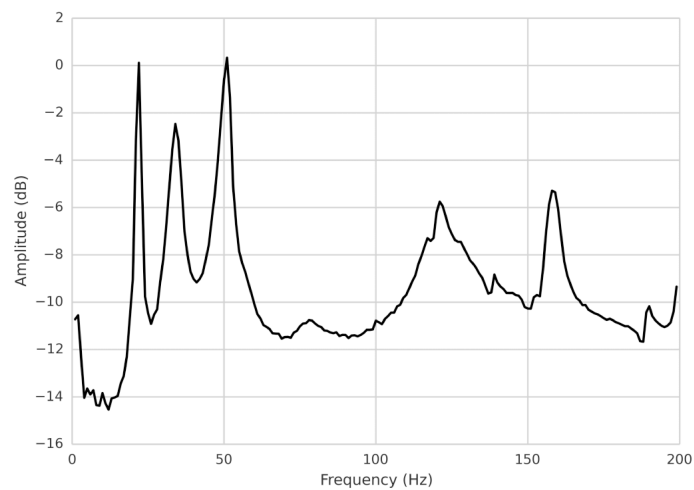

Figure 3 – The average frequency response of the transformer

### 3.2 FEA Results and Discussion

The Abaqus software was used for FEA simulation, with the preliminary experimental result as a reference. The anisotropic property is set to the silicon steel according to previous research (2). Table 1 presents the first four natural frequencies obtained from the FEA simulation and the experimental measurements. Figure 4 presents the mode shape measured from experiment and simulated by FEA. As shown in Table 1 and Figure 4, the FEA result is close to the measurement result in both from the natural frequencies and the mode shape. Therefore, this FEA model will be used in later FEA-FRFs calculations.

Table 1 – The natural frequency of the transformer

| Mode | Measured Frequency (Hz) | FEA Frequency (Hz) | Relative Error |
|------|------------------------|--------------------|----------------|
| 1    | 23                     | 24.5               | 6.5%           |
| 2    | 35                     | 32.8               | 6.3%           |
| 3    | 50                     | 51.5               | 6.5%           |
| 4    | 129                    | 128.9              | 0.07%          |

After establishing the FEA model, the next step is to calculate the FEA-FRF at the seed frequency, which is 100 Hz in this case. In this FEA model, $N = 635$ nodes of the core are selected as the force nodes and response nodes and each of them has three directions. Thus, there are a totally 3,629,025 ($9N^2$) FEA-FRFs that need to be calculated. A total of 8 hours were spent on producing this dataset on a computer with an Intel i7 processor with 8GB of RAM.

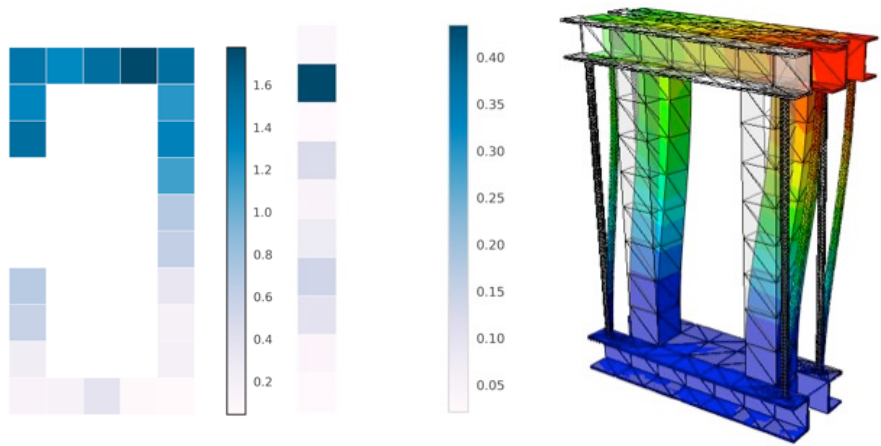

Figure 4 – The experimental response and the corresponding FEA mode shape at 23 Hz

## 4. EXTRATREE AND THE GA

After obtaining the full dataset of FEA-FRFs, the GA and ExtraTree were applied to them in order to obtain the best Green's function model. In this section, the procedure for selecting the ExtraTree algorithm are presented and the details of parameters and rules for the GA are given.

### 4.1 The Selection of ExtraTree

Before the GA is executed, a regression model need to be selected for the fitness calculation. The dataset was trained and tested to observe the performances of traditional machine learning algorithms for the best one to be selected. In this case, the dataset was separated into a training set and a testing set with the same sample sizes. Each traditional machine learning algorithm was trained on the training set and tested on the testing set. Their performances are presented in Table 2, and include the correlation coefficient ($R$), the relative root-mean-square error (RMSE-r) and the time consumed. To prevent from over fitting, all the algorithms were trained without tuning any parameters. The features used in training include the coordinates of the force node and the response node ($x, y, z$) and the directions of the force and response (1-3). The importance values of the features are shown in Table 3.

As shown in Table 2, ExtraTrees performed best with the highest $R$ and the lowest RMSE-r, with acceptable calculation time consumed. Therefore, ExtraTree is selected among traditional machine learning algorithms that were used in subset selection and development of the final Green's function model development.

Table 2 – The performances and calculation times of the models

| Model | R | RMSE-r | Time (s) |
|---|---|---|---|
| Linear | 0.212 | 0.964 | 5.94 |
| Bayesian Ridge | 0.212 | 0.964 | 8.23 |
| Random Forest | 0.99981 | 0.0190 | 11.32 |
| K-Nearest Neighbourhood | 0.997 | 0.0787 | 2299.5 |
| ExtraTree | 0.9999 | 0.0138 | 10.88 |
| AdaBoost | 0.978 | 0.398 | 70.55 |

Table 3 – The importance values of the features of the ExtraTree model

| Feature | Force Direction | Response Direction | $x_F$ | $y_F$ | $z_F$ | $x_R$ | $y_R$ | $z_R$ |
|---|---|---|---|---|---|---|---|---|
| Importance | 0.198 | 0.143 | 0.170 | 0.116 | 0.019 | 0.165 | 0.167 | 0.023 |

## 4.2 The Implementation of the GA

In the GA, a few rules and parameters must be set. These include the format of the individual, the fitness criterion, the crossover method, and the parameters such as the population, the remaining ratio, the random selection ratio, and the mutation possibility.

All the parameters were selected using the convergence test, which compared the convergence speed of the GA for minimising the FEA-FRF quantity. To reflect the FRFs selection of the FRFs, the individual is designed as the selection matrix $S$ with dimensions 635 * 635. Each element inside $S$ represents the selection status of the corresponding FRF. When the value $S_{ij}$ equals to 1, it means the FRF from node $i$ to node $j$ is selected in the ExtraTree model training and if it is equal to 0, the corresponding FRF is not selected.

For the crossover method, the *and* operation is chosen. The father and mother individuals are selected randomly from the parent set to generate the next generation, as shown in Figure 5.

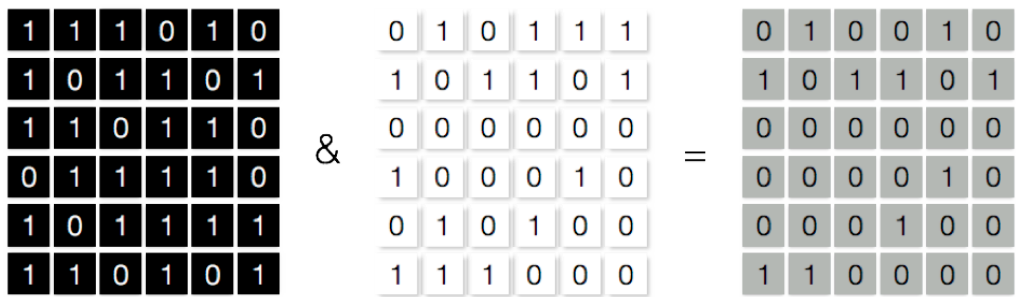

Figure 5 – The *and* crossover method

To define the fitness, variate methods were tried. The best fitness setting is presented in Eq. 2.

$$f = \frac{1}{rq} H[r_t - r] \quad (2)$$

where $r$ is the RMSE-r of the model, $r_t$ is the maximum tolerated $r$, $q$ represents the selected FEA-FRF quantity, and $H$ is the Heaviside step function. Different FEA-FRFs subsets will be obtained by setting different values of $r_t$ in the GA at the seed frequency.

Finally, based on the results from the convergence test, the remaining ratio, the random selection ratio and the mutation possibility are 0.4, 0.1, and 0.02, respectively.

After defining the rules and parameters, the GA is carried out on the full dataset of the FEA- FRF to obtain the best subset. To calculate the Green's function at other frequencies more efficiently with less calculation time, the subset FEA-FRFs at other frequencies are calculated and the results inputted to the ExtraTree algorithm for training. The trained model is the Green's function model at the corresponding frequencies.

## 5. RESULT AND DISCUSSION

In this paper, the seed frequency is 100 Hz. The $r_t$ values are set to be 0.8, 0.7, 0.5, and 0.3 dB. Their corresponding FEA-FRF calculation reduction times, which is the whole FEA-FRF quantity over the quantity of the selected subset, come out to be 5.0, 4.3, 2.2, and 2.0, respectively. Their corresponding optimised FEA-FRF subsets are applied at the 50 Hz, 70 Hz, 98 Hz, 102 Hz, 120 Hz, and 150 Hz, which cover the natural frequencies and others. With the model developed on these frequencies, their RMSE-r values on the full dataset are calculated and presented in the left panel of Figure 6. To show the GA selected FRFs are a better combination than the randomly selected one, the models based on random FRF subset are compared and the corresponding RMSE-r values are presented in the right panel of Figure 6. The random FRFs' subsets have the same quantities as the ones they are comparing.

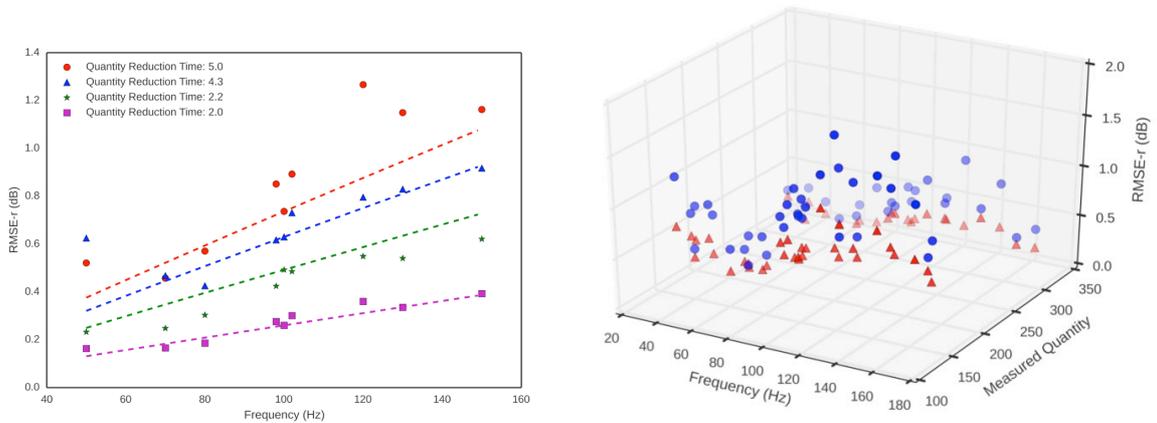

Figure 6 – The RMSE-r values of the FRF trained model

In the right panel of Figure 6, the red dots are the RMSE-r values of the models built on the GA selected subsets and the blue dots are those trained on randomly selected subsets. It is clear that all the models built on the GA selected subset have smaller RMSE-r values than those on the randomly selected subset, with an average 54.5% reduction. It is also obvious that, with the decrease of the calculation time, the models' error will increase. Another important phenomenon is that the RMSE-r values decrease for the models whose frequency is lower than the GA seed frequency, while they increase for the model where it is higher, which indicates the relationship between the error and frequency.

To explain this phenomenon, a domain zone theory is proposed in this paper. Each frequency has an associated domain zone such that the variation of this zone can be represented by a selected node inside it. With increasing frequency, the size of the zone will decrease as the wavelength becomes smaller. A smaller zone means more zones in the transformer core, which will require more nodes to represent. Therefore, the GA selected FEA-FRFs cover the zones at the frequencies of 100 Hz and below. In other words, the GA can identify the domain zones under the seed frequency and represented by the FEA-FRF subset. Therefore, the GA selected FRF subset can be applied to the seed frequency and below. When used to estimate the values at higher frequencies, a larger error will be produced. The left panel of Figure 7 presents the actual model errors. Based on an understanding of the domain zone theory, the estimated error is given by:

$$r_e = r_s \frac{f_c}{f_e} \qquad (3)$$

Where $f_c$ is the calculated frequency, $f_s$ is the seed frequency, $r_e$ is the RMSE-r of the $f_c$, and $r_s$ is the RMSE-r on the $f_s$. The estimated error is indicated in the left panel of Figure 6 by the dashed lines. It is clear from Eq. 3, the error is proportional to the frequency, and this is verified by the high correlation between the dots and dashed line in Figure 6. This agrees with the domain zone assumption.

Therefore, this method was demonstrated to be useful in speeding up the FEA calculation for a Green's function with controlled error. The ExtraTree Green's function model can be trained to predicted any FEA-FRF. The calculation time can be reduced by reducing the FEA-FRFs to be calculated in FEA and the relative error can be estimated according to Eq. 3. This method can be applied not only to a transformer core but also to other complex structures, which will be studied in later research. In addition, the subset from the GA will be measured on an actual transformer for a more precise Green's function model, which will be introduced in a future paper.

## 6. CONCLUSION

In this paper, a novel method for Green's function modelling was proposed by combining the FEA and the ExtraTree algorithm. It was shown that the Green's function can be modelled by an ExtraTree model trained on a subset of FEA-FRFs. A subset of the total number of FEA FRFs was obtained by a genetic algorithm at a seed frequency. This can be applied to other frequencies so that the time required to calculate the Green's function of a structure can be efficiently reduced through calculate fewer finite element and additional ExtraTree training with a controlled error. Further work can be performed by combining this method with experiment to increase the model accuracy.

## 7. ACOKNOWLEDGEMENT

The authors would like to show their gratitude for the support from the International Postgraduate Research Scholarship and the Nectar cloud computation service in Australia and support from Zhejiang University. The authors also acknowledge the NSFC (Grant Nos. 11574269) for support.